# LINGUISTIC-MATHEMATICAL STATISTICS IN REBUS, LYRICS, JURIDICAL TEXTS, FANCIES AND PARADOXES

Florentin Smarandache
University of New Mexico
200 College Road
Gallup, NM 87301, USA
E-mail: smarand@unm.edu

**Abstract.**
This is a collection of linguistic-mathematical approaches to Romanian rebus, poetical and juridical texts, and proposes fancies, recreational math problems, and paradoxes. We study the frequencies of letters, syllables, vowels in various poetry, grill definitions in rebus, and rebus rules. We also compare the scientific language, lyrical language, and puzzles' language, and compute the Shannon entropy and Onicescu informational energy.

## INTRODUCTION

The aim of this section is the investigation of some combinatorial aspects of written language, within the framework determined by the well-known game of crossword puzzles. Various types of probabilistic regularities appearing in such puzzles reveal some hidden, not well-known restrictions operating in the field of natural languages. Most of the restrictions of this type are similar in each natural language. Our direct concern will be the Romanian language.

Our research may have some relevance for the phono-statistics of Romanian. The distribution of phonemes and letters is established for a corpus of a deviant morphological structure with respect to the standard language. Another aspect of our research may be related to the so-called tabular reading in poetry. The correlation horizontal-vertical considered in the first part of the paper offers some suggestions concerning a bi-dimensional investigation of the poetic sing.

Our investigation is concerned with the Romanian crossword puzzles published in [4]. Various concepts concerning crossword puzzles are borrowed from N. Andrei [3]. Mathematical linguistic concepts are borrowed from S. Marcus [1], and S. Marcus, E. Nicolau, S. Stati [2].

## SECTION 1. THE GRID

### §1. MATHEMATICAL RESEARCHES ON GRIDS

It is known that a word in a grid is limited on the left and right side either by a black point or by a grid final border.

We will take into account the words consisting of one letter (though they are not clued in the Rebus), and those of two (even they have no sense (e.g. N T, RU,…)), three or more letters – even they represent that category of rare words (foreign localities, rivers, etc., abbreviations, etc., which are not found in the Romanian Language Dictionary (see [3], pp. 82-307 ("Rebus glossary"))).

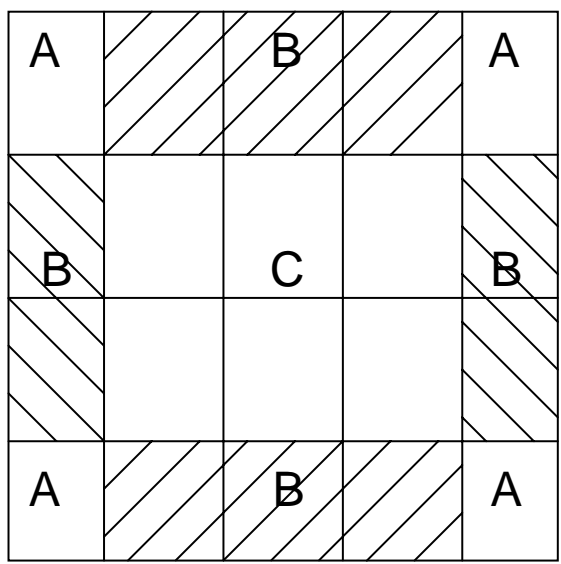

The grids have both across and down words.
We divide the grid into 3 zones:
a) the four peaks of the grid (zone A)
b) grid border (without de four peaks) (zone B)
c) grid middle zone (zone C)
We assume that the grid has $n$ lines, $m$ columns, and $p$ black points.
Then:

**Proposition 1.** The words overall number (across and down) of the grid is equal to $n + m + pNB + 2 \cdot pNC$, where

$pNB$ = black points number in zone $B$,

$pNC$ = black points number in zone $C$.

*Proof:* We consider initially the grid without any black points. Then it has $n + m$ words.

- If we put a black point in zone $A$, the words number is the same. (So it does not matter how many black points are found in zone A).

- If we put a black point in zone $B$, e.g. on line 1 and column $j$, $i < j < m$, words number increases with one unit (because on line 1, two words were formed (before there was only one), and on column $j$ one word rests, too). The case is analog if we put a black point on column 1 and line $i$, $1 < i < n$ (the grid may be reversed: the horizontal line becomes the vertical line and vice versa). Then, for each point in zone B a word is added to the grid words overall number.

- If we put a black point in zone $C$, let us say $i$, $1 < i < n$, and column $j$, $1 < j < m$, then the words number increases by two: both on line $i$ and column $j$ two words appear now, different from the previous case, when only one word was there on each line. Thus, for each black point in zone $C$, two words are added at the grid words overall number. From this proof results:

**Corollary 1**. Minimum number of words of grid $n \times m$ is $n + m$. Actually, this statement is achieved when we do not have any black points in zones $B$ and $C$.

**Corollary 2.** Maximum number of words of a grid $n \times m$ having $p$ black points is $n + m + 2p$ and it is achieved when all $p$ black points are found in zone $C$.

**Corollary 3.** A grid $n \times m$ having $p$ black points will have a minimum number of words when we fix first the black points in zone $A$, then in zone $B$ (alternatively – because it is not allowed to have two or more black points juxtaposed), and the rest in zone $C$.

**Proposition 2**. The difference between the number of words on the horizontal and on the vertical of a grid $n \times m$ is $n - m + pNBO - pNBV$, where

$pNBO$ = black points number in zone $BO$,

$pNBV$ = black points number in zone $BV$.

We divide zone $B$ into two parts:
- zone $BO$ = $B$ zone horizontal part (line 1 and $n$)
- zone $BV$ = $B$ zone vertical part (line 1 and $m$).

The proof of this proposition follows the previous one and uses its results.

If we do not have any black points in the grid, the difference between the words on the horizontal and those on the vertical line is $n - m$.

- If we have a black point in zone $A$, the difference does not change. The same for zone $C$.

If we have a black point in zone $BO$, then the difference will be $n - m - 1$. From this proposition 2 results:

**Proposition 3**. A grid $n \times m$ has $n + pNBO + pNC$ words on the horizontal and $m + pNBV + pNC$ words on the vertical.

The first solving method uses the results of propositions 1 and 2.

The second method straightly calculates from propositions 1 and 2 the across and down words number (their sum (proposition 1) and difference (proposition 2) are known).

**Proposition 4.** Words mean length (=letters number) of a grid $n \times m$ with $p$ black points is $\geq \dfrac{2(nm - p)}{n + m + 2p}$.

Actually, the maximum words number is $n + m + 2p$, the letter number is $nm - p$, and each letter is included in two words: one across and another down. One grid is the more crossed, the smaller the number of the words consisting of one or two letters and of black points (assuming that it meets the other known restrictions). Because in the Romanian grids the black points percentage is max.

15% out of the total (rounding off the value at the closer integer – e.g. 15% with a grid 13x13 equals 25.35 $\approx$ 25; with a grid 12x12 is 21.6 $\approx$ 22), so for the previous properties, for grids $n \times m$ with $p$ black points we replace $p$ by $\left[\dfrac{3}{20}\right]nm$, where

$[x] = \max\{\alpha \in \mathbb{N}, \ |\alpha - x| \leq 0.5\}$.

## §2. STATISTIC RESEARCHES ON GRIDS

In [1] we find the notion "écart of a sound x", denoted by $\alpha(x)$, which equals the difference between the rank of $x$ in Romanian and the rank of $x$ in the analyzed text.

We will extend this notion to the notion of *a text écart* which will be denoted by: $\alpha(t)$, and

$$\alpha(t) = \frac{1}{n}\sum_{i=1}^{n}|\alpha(A_i)|$$

where $\alpha(A_i)$ is $A_i$ sound écart (in [1]) and $n$ represents distinct sounds number in text $t$. (If there are letters in the alphabet, which are not found in the analyzed text, these will be written in the frequency table giving them the biggest order.)

**Proposition 1.** We have a double inequality:

$$0 \leq \alpha(t) \leq \frac{n-1}{2} + \frac{1}{n}\left[\frac{n}{2}\right]$$ where $[y]$ represents the whole part of real number $y$.

Actually, the first inequality is evident.

Let $\Phi = \begin{pmatrix} 1 & 2 & \ldots & n \\ j_1 & j_2 & \ldots & j_n \end{pmatrix}$. Then $\sum_{i=1}^{n}|\alpha(A_i)| = \sum_{i=1}^{n}|i - j_i|$

This permutation constitutes a mathematical pattern of the two frequency tables of sounds; in Romanian (the first line), in text t (the second line).

For permutation $\psi = \begin{pmatrix} 1 & 2 & \ldots & n-1 & n \\ n & n-1 & \ldots & 2 & 1 \end{pmatrix}$ we have

$$\sum_{i=1}^{n}|i - j_i| = 2[(n-1) + (n-3) + (n-5) + \ldots] = 2\sum_{k=1}^{\left[\frac{n}{2}\right]}(n - 2k + 1) =$$

$$= 2\left[\frac{n}{2}\right]\left(n - \left[\frac{n}{2}\right]\right) = \frac{n(n-1)}{2} + \left[\frac{n}{2}\right],$$

where $\alpha(t) = \frac{n-1}{2} + \frac{1}{n} \cdot \left[\frac{n}{2}\right]$.

By induction with respect to $n \geq 2$, we prove now the sum $S = \sum_{i=1}^{n}|i - j_i|$ has max. value for permutation $\psi$.

For $n = 2$ and 3 it is easily checked directly. Let us suppose the assertion true for values $< n + 2$. Let us show for $n + 2$:

$$\psi = \begin{pmatrix} 1 & 2 & \ldots & n+1 & n+2 \\ n+2 & n+1 & \ldots & 2 & 1 \end{pmatrix}$$

Removing the first and last column, we obtain:

$$\psi' = \begin{pmatrix} 2 & \ldots & n+1 \\ n+1 & \ldots & 2 \end{pmatrix},$$

which is a permutation of $n$ elements and for which $S$ will have the same value as for permutation

$$\psi'' = \begin{pmatrix} 1 & \dots & n \\ n & \dots & 1 \end{pmatrix},$$

i.e. max. value ($\psi''$ was obtained from $\psi'$ by diminishing each element by one).

The permutation of 2 elements $\eta = \begin{pmatrix} 1 & n+2 \\ n+2 & 1 \end{pmatrix}$ gives maximum value for $S$.

But $\psi$ is obtained from $\psi'$ and $\eta$;

$$\psi(i) = \begin{cases} \psi'(i), & \text{if } i \notin \{1, n+2\} \\ \eta(i), & \text{otherwise} \end{cases}$$

*Remark* : The bigger one text écart, the bigger the "angle of deviation" from the usual language.

It would be interesting to calculate, for example, the écart of a poem.

Then the notion of écart could be extended even more:

a) *the écart of a word* being equal to the difference between word order in language and word order in the text;
b) *the écart of a text (ref. words):*

$$\alpha_c(t) = \frac{1}{n} \sum_{i=1}^{n} |\alpha_c(a_i)|,$$

where $\alpha_c(a_i)$ is word $a_i$ écart, and $n$ - distinct words number in the text $t$.

*

We give below some rebus statistic data. By examining 150 grids [4] we obtain the following results:

*Occurrence frequency of words in the grid, depending on their length (in letters)*

| Letter order | Letter | Letter occurrence mean percentage | Vowels mean percentage | Consonants mean percentage |
|---|---|---|---|---|
| 1 | A | 15.741% | | |
| 2 | I | 12.849% | | |
| 3 | T | 9.731% | | |
| 4 | R | 9.411% | | |
| 5 | E | 8.981% | | |
| 6 | O | 5.537% | | |
| 7 | N | 5.053% | | |
| 8 | U | 4.354% | 47.462% | 52.538% |
| 9 | S | 4.352% | | |
| 10 | C | 4.249% | | |
| 11 | L | 4.248% | | |
| 12 | M | 4.010% | | |
| 13 | P | 3.689% | | |
| 14 | D | 1.723% | | |
| 15 | B | 1.344% | | |
| 16 | G | 1.290% | | |
| 17 | F | 0.860% | | |
| 18 | V | 0.806% | | |
| 19 | Z | 0.752% | | |
| 20 | H | 0.537% | | |
| 21 | X | 0.430% | | |
| 22 | J | 0.053% | | |
| 23 | K | 0.000% | | |

It is easy to see that a percentage of 49,035% consists of the words formed only of 1, 2 or 3 letters; - of course, there are lots of incomplete words.

*

The study of 50 grids resulted in:
  *Occurrence frequency of words in a grid* (see next page).
It is noticed that vowels percentage in the grid (47.462%) exceeds the vowels percentage in language (42.7%).
So, we can generalize the following:
  *Statistical proposition* (1): In a grid, the vowels number tends to be almost equal to 47.5% of the total number of the letters.
  Here is some evidence: one word with $n$ syllables has at least $n$ vowels (in Romanian there is no syllable without vowel (see [2]).
  The vowels percentage in Romanian is 42.7%; because a grid is assumed to form words across and down, the vowels number will increase. Also, the last two lines and columns are endings of other words in the grid; thus they will usually have more vowels. When black points number decreases, vowels number will increase (in order to have an easier crossing, you need either more black points or more vowels) (A vowel has a bigger probability to enter in the contents of a word than a consonant.)

Especially in "record grids" (see [3], pp. 33-48) the vowels and consonants alternation is noticed. Another criterion for estimating the grid value is the bigger deviation from this "statistical law" (the exception confirms the rule!): i.e. the smaller the vowel percentage in a grid, the bigger its value.

*Statistical proposition* (2): Generally, the horizontal words number 73 equals the vertical one.

Here is the following evidence: 100 classical grids were experimentally analyzed, in [4], getting the percentage of 49.932% horizontal words. Usually, the classical grids are square clues, the difference between the horizontal and vertical words being (see Proposition 2):
$$n - m + pNBO - pNBV = pNBO - pNBV.$$

The difference between the black points number in zone $BO$ and zone $BV$ can not be too big ($\pm 1$, $\pm 2$ and rarely $\pm 3$). (Usually, there are not many black points in zone B, because it is not economical in crossing (see proof of Proposition 1)).

Taking from [1] the following letters frequency in language:

| | | | | | |
|---|---|---|---|---|---|
| 1. E | 5. N | 9. L | 13. P | 17. G | 21. J |
| 2. I | 6. T | 10. S | 14. M | 18. F | 22. X |
| 3. A | 7. U | 11. O | 15. B | 19. Z | 23. K |
| 4. R | 8. C | 12. D | 16. V | 20. H | |

(because in the grid Ă, Â, Î, Ş, Ţ: are replaced by A: I: S: T, respectively, in the above order they were cancelled) the écart of the 150 grids becomes
$$\alpha(g) = \frac{1}{23} \sum_{i=1}^{23} |\alpha(A_i)| \approx 1.391;$$
the entropy is:
$$H_1 = -\frac{1}{\log_{10} 2} \sum_{i=1}^{23} p_i \log_{10} p_i \approx 3.865$$
and the informational energy (after O. Onicescu) is:
$$E(g) = \sum_{i=1}^{23} p_i^2 \approx 0.084$$
Examining 50 grids we obtain:

*Words frequency in a grid with respect to the syllables number*

| Mean percentage of occurrence of a word in a grid | | | | | | | | Mean length of a word in syllables |
|---|---|---|---|---|---|---|---|---|
| 1 syllable | 2 | 3 | 4 | 5 | 6 | 7 | 8 | |
| 35.588% | 26.920% | 21.765% | 9.551% | 5.294% | 0.882% | 0.000% | 0.000% | 2.246 |

(in the category of the one syllable-words, the word of one, two or, three letters, without any sense – rare words – were also considered.) One can see that the percentage of words consisting of one and two syllables is 65.508% (high enough).

Another statistics (of 50 grids), concerning the predominant parts of speech in a grid has established the following first three places:
1. nouns 45.441%
2. verbs 6.029%
3. adjectives 2.352%

Notice the large number of nouns.

\*

**SECTION II. REBUS CLUES**

**§1. STATISTICAL RESEARCHES ON REBUS CLUES**

Studying the clues of 100 "clues grids", the following statistical data resulted:
*Rebus clues frequency according to their length (words number)*
(see the next page)
It is noticed that the predominant clues are formed of 2, 3, or 4 words. For results obtained by investigating 100 "clues grids", see the next page.

It is worth mentioning that vowels percentage (46.467%) from rebus clues exceeds vowels percentage in the language (42.7%).

By calculating the clues écart (in accordance with the previous formula) it results:

$$\alpha(dr) = \frac{1}{27}\sum_{i=1}^{27}|\alpha(A_i)| \approx 1.185$$

(sound frequency used by Solomon Marcus in [1] was used here), the entropy (Shannon) is:

$$H_1 = -\frac{1}{\log_{10} 2} \sum_{i=1}^{27} p_i \log_{10} p_i \approx 4.226$$

and informational energy (O. Onicescu) is:

$$E(dr) = \sum_{i=1}^{27} p_i^2 \approx 0.062.$$

(The calculations were done by means of a pocket calculator ).

*Letters occurrence frequency in the rebus clues*

| Letter order | Letter | Mean percentage of letter occurrence in clues | Vowels percentage | Consonants mean percentage | Letters no. (mean) necessary to clue a grid | Mean length of a word (in letters) used in clues |
|---|---|---|---|---|---|---|
| 1 | E | 10.996% | | | | |
| 2 | I | 9.778% | | | | |
| 3 | A | 9.266% | 46.679% | 53.321% | 657.342 | 4.374 |
| 4 | R | 7.818% | | | | |
| 5 | U | 6.267% | | | | |
| 6 | N | 6.067% | | | | |
| 7 | T | 5.611% | | | | |
| 8 | C | 5.374% | | | | |
| 9 | L | 4.920% | | | | |
| 10 | O | 4.579% | | | | |
| 11 | P | 4.027% | | | | |
| 12 | Ă | 3.992% | | | | |
| 13 | S | 3.831% | | | | |
| 14 | Î | 3.309% | | | | |
| 15 | D | 3.079% | | | | |
| 16 | Â | 1.801% | | | | |
| 17 | V | 1.527% | | | | |
| 18 | F | 1.449% | | | | |
| 19 | Ş | 1.360% | | | | |
| 20 | Ţ | 1.338% | | | | |
| 21 | G | 1.330% | | | | |
| 22 | B | 1.238% | | | | |
| 23 | H | 0.532% | | | | |
| 24 | J | 0.358% | | | | |
| 25 | Z | 0.092% | | | | |
| 26 | X | 0.037% | | | | |
| 27 | K | 0.024% | | | | |

The Craiova University
Natural Sciences Department

# SECTION III. HYPOTHESIS ON THE DETERMINATION OF A RULE FOR THE CROSS WORDS PUZZLES

The problems of cross words are composed, as we know, of grids and definitions. In the Romanian language one imposes the condition that the percentage of black boxes compared to the total number of boxes of the grid not to go over 15%.

Why 15%, and not more or less? This is the question to which this article tries to answer. (This question is due to Professor Solomon MARCUS - National Symposium of Mathematiques "Traian Lalesco", Craiova University, June 10, 1982).

First of all we present here a table which shows in a synthetic manner, a statistics on the grids containing a very small percentage of black boxes (of [2], pp. 27-29):

## THE GRIDS-RECORDS

| Grid dimension | Minimum number of registered black boxes | Percentage of black boxes | Number of grids-records constructed until June 1, 1982 |
|---|---|---|---|
| 8x8 | 0 | 0.000% | 24 |
| 9x9 | 0 | 0.000% | 3 |
| 10x10 | 3 | 3.000% | 2 |
| 11x11 | 4 | 3.305% | 1 |
| 12x12 | 8 | 5.555% | 1 |
| 13x13 | 12 | 7.100% | 1 |
| 14x14 | 14 | 7.142% | 1 |
| 15x15 | 17 | 7.555% | 1 |
| 16x16 | 20 | 7.812% | 2 |

In this table, one can see that the larger the dimension of the grid, the larger is the percentage of black boxes, because the number of long words is reduced.

The current dimensions for grids go from 10x10 to 15x15.

One can notice that the number of the grids having a percentage of black boxes smaller than 8 is very reduced: the totals in the last column represent all the grids created in Romania since 1925 (the appearance of the first problems of cross words in Romania), until today. It is thus seen that the number of the grid-records is negligible when one compares it with the thousands of grids created. For this reason, the rule that imposed the percentage of the black boxes, should have established to be greater than 8%. But the cross words being puzzles, they must address to a large audience, thus one did not have to make these problems too difficult.

From which a percentage of black boxes at least equal to 10%.

They must be not too easy either, that is not to necessitate any effort from those who would compose them, from where a percentage of black boxes smaller than 20%. (If not, in effect, it becomes possible to compose grids wholly formed of words boxes of 2 or 3 letters).

To support the second assertion, one assumes that the average length of the words of a $n \times m$ grid with $p$ black boxes is sensible equal to $\dfrac{2(n \cdot m - p)}{n + m + 2p}$ (from [3]. § 1, Prop. 4). For us, $p$ is 20% of $n \cdot m$, therefore it results that

$$\frac{2(n \cdot m - \dfrac{20}{100} n \cdot m)}{n + m + 2 \dfrac{20}{100} n \cdot m} \leq 3 \Leftrightarrow \frac{1}{n} + \frac{1}{m} \geq \frac{2}{15}.$$

Thus, for current grids having 20% of black boxes, the average lengths of the words would be smaller than 3.

Similarly at the beginnings of the puzzle of cross words the percentage of black boxes were not too large: thus in a grid from 1925 of 11x11, one counts 33 black boxes, therefore a percentage of 27.272% (from [2], p. 27).

While being developed, for these puzzles were imposed "stronger" conditions – that is a reduction in the black boxes.

For selecting a percentage between 10 and 20%, it is supposed that the peoples' predilection for round numbers was essential (the cross words are puzzles, no need for mathematic precision of sciences). That's why the rule of 15%.

A statistic (from [3], § 2), shows that the percentage of black boxes in the current grids is approximately 13.591%. The rule is thus relatively easy to follow and it can only attract new crossword enthusiasts.

To completely answer the proposed question, one would need to consider also some philosophical, psychological, and especially sociological aspects, especially those connected to the history of this puzzle, its ulterior development, and with its traditions.

SECTION IV. THE LANGUAGE OF SPIRITUAL REBUS DEFINITIONS

"The rebus' language" is somewhere at the border of the scientific language and, that, perhaps, having many common things with usual language too, and even with the musical one (the puzzles, because they have a certain acoustic resonance).

While the semantic deficiencies, having direct definitions (close to those from dictionary [3], pp. 50-56) of a language close to the scientific one (even to the usual one through the simple mode of expression) of "the grid's definitions". The language is close to the poetic one. There are even literary definitions (see [3], p. 57, [4]), which utilize literary stylistic procedures: like the metaphor, the comparison, the allegory, practice, etc. Later we will present a parallelism between the SCIENTIFIC LANGUAGE, POETIC LANGUAGE, REBUS' LANGUAGE ("THE GRIDS' DEFINITIONS") closely following the rules from [1] (chap. "Oppositions between the scientific language and the poetic one"), results which we will limit to the rebus' language.

| SCIENTIFIC LANGUAGE | POETIC LANGUAGE | REBUS' LANGUAGE |
|---|---|---|
| - rational hypothesis | - emotional hypothesis | - rational + emotional hypothesis (reading the definition, you think for an instant, sometimes you go on a wrong road; when you err the answer (the corresponding word from the grid, you get enlightened and enthusiast). |
| - logical density | - density of suggestion | - logical density + suggestion (the definition must use very few words to explain a lot – logical density); to be unpublished, enlightening, emotional (density of suggestion). |
| - infinite synonymy | - absent synonymy | - reduced synonymy (not truly infinite, but not absurd); (two identical words from the grid cannot have more than one rebus definition; but a definition will be almost uniquely expressed, therefore the synonymy is quasi absent). |
| - absent anonymity | - infinite anonymity | - large anonymity (neither absent, nor infinite) (in the case of the definition, the meaning is up to the author: |
|   |   | even if the reader understands something else, it will intervene the rational part, the word must fulfill the proper place in the grid, even the literary definitions, in the grids, don't have anymore an infinite anonymity, because here intervene also the rational part: the finding by all means of an answer: in the case of the theme grids with direct definitions, the anonymity is almost absent). |
| - artificial | - natural | - natural and artificial (in general the definitions have a natural character; but the definitions based on letter's puzzles (example, the definition "Night's beginning" has the answer "NI" have an artificial character). |
| - general | - singular | - singular and general (only the definitions based on the puzzles of letters may have a general character). |
| - translatable | - untranslatable | - translatable (in the sense that the definition has a logical meaning). |
| - the presence of style problems | - the absence of style problems | - the absence of style problems (the same definition cannot be used without changing the nuance – while a word in the grid can be defined in multiple ways). |
| - finitude in space, constant in time | - variability in space and time | - the variability in space and time, smaller variability than that from the poetic language. |
| - numerable | - innumerable | - innumerable |
| - transparent | - opaque | - semi-opaque (or semitransparent - at the |

|  |  | beginning the definition seems opaque, until one finds the answer). |
| --- | --- | --- |
| - transitive | - reflexive | - reflexive (except, again, the definitions based on games of letters, which have also a transitional character). |
| - independency on expression | - dependency on expression | - dependency on expression. |
| - independency on musical structure | - dependency on musical structure | - dependency on musical structure. |
| - paradigmatic | - syntagmatic | - syntagmatic |
| - concordance between the paradigmatic and syntagmatic distance | - non concordance between the paradigmatic and syntagmatic distance | - the paradigmatic and syntagmatic distance (are pairs of different words, word games, methods used ass in poetry). |
| - short contexts | - long contexts | - short contexts (1) (here it is closer to the scientific language, because it is taken into account the Latin proverb "*Non multa sed multum*"; from the anterior statistic investigations it resulted that the medium length of a (spiritual) rebus definition is 4.192 words: the definitions with letter puzzles usually have very few words. |
| - contextual dependency | - it tends towards expression independency | - contextual dependency (in the case of the theme grids it is also a small dependency; there exist also rare cases when a definition is dependent of an anterior definition (usually the definitions with letters or word games)). |
| - logic | - illogic | - logic |
| - denotation | - annotation | - connotation (if a definition would reveal the direct meaning of an word, we would have direct definitions (like in a |
|  |  | dictionary)) and then we would totally loose "the surprise", "the spirituality", "the ingenious", "the spontaneity" of thematic grids, the definitions with denotative character. |
| - routine | - creation | - creation and … experience (not to call it routine!) |
| -general stereotypes | - personal stereotypes | - personal stereotypes (it exists even the so called grids of "personal manner" – (see [3], pp. 56-58) |
| - explicable | - ineffable | - ineffable … which explains it! (Taken separately, the definition, not-seen as a question, is ineffable taken along, with the answer becomes explicable: in general, the definition presents also an ambiguity degree (more tracks for guidance) – otherwise it would be banal – a degree of indetermination: it is used many times the proper sense instead of the figurative one, or reciprocally defined it has also its own logic, which becomes tangible once one finds the answer). |
| - lucidity | - magic | - magic – lucidity (in accordance with those that are immediately anterior) (at the beginning the rebus language dominates the person, until he finds the "key" when he'll become at his turn the dominant – the poetic language. |
| - predictable | - unpredictable | - at the beginning is unpredictable, and becomes predictable after solving it: (unpredictable converted in predictable) . |

# CONSIDERATIONS REGARDING THE SCIENTIFIC LANGUAGE AND "LITERARY LANGUAGE"

As in nature nothing is absolute, evidently there will not exist a precise border between the scientific language and "the literary" one (the language used in literature): thus there will be zones where these two languages intersect.

In [1], chapter "Instances between the scientific and poetic languages", Solomon Marcus presents the differences between these two, differences that make them closer.

We will skate a little on the edge of this material, presenting common parts of the scientific language and the literary language:
- both are geared to find the unpublished, the novelty
- both suppose a creative process (finding the solution of a problem means creation: writing of a phrase the same).
- both literature and science have an art of being taught, studied and learned (the methodology of teaching arithmetic, or Romanian language, etc.) .
- in science too there is an esthetic (for example: "the mathematical esthetic"), the same in literature there exists a logic (even the absurd of Eugene Ionesco, the myths of Mircea Eliade have their own specific logic: analogously, we can extend the idea to Tristan Tzara's Dadaism, which has a specific logic (of construction; one cuts words from newspapers, mix them, and then form verses).
- the scientific development implies a literary development in a special sense: it appeared, thus, the science-fiction literature in literary writings which use informations obtained by science: contemporaneous literature treats also scientific problems (for example Augustin Buzura wrote the roman "The absents" describing the life of a medical researcher: the engineer poet George Stanca introduces technical terms in his poems; one verse from his volume "Maximum tenderness" sounds: "$\sin^2 x + \cos^2 x = 1$"!); analogously the engineer poet Gabriel Chifu (the volume "An interpretation of the Purgatory") and mathematics professor Ovidiu Florentin, author of a volume even entitled "Formulas for the spirit" – each poem being considered as a momentous "formula" (depending of time, place, space, individual) for the spirit.
- even the writing of some contemporary novels inspired from the worker's and peasant's life requires a scientific documentation from the writers' part.

The literature has an esthetic influence for science; there exist mathematical metaphors (see [1], [2]) and, in general, we can say "scientific metaphors", one cannot know what ideas and relations will be discovered in science. The understanding degree (exegesis) of a poetry and of a literary text in general, depends also of the culture's degree of each individual, of his initiation (the seniority in that domain), of his scientific knowledge.
- there are many scientists who, besides their scientific works, write also literary works or related domains (for example, the memories book of the academician (mathematician) Octav Onicescu "On the life's roads", the renown Romanian physician Gheorghe Marinescu writes poems (using Dacic words), under the penname George Dinizvor, the great Ion Barbu – Dan Barbilian excelled as a poet and as a mathematician. The great poet Vasile Voiculescu was a good physician; and the mathematics professor Aurel M. Buricea writes poetry, analogously the mathematician Ovidiu Florentin –

Florentin Smarandache writes poems and mathematics articles; in the world literature we find the poet-mathematician Omar Khayyam and Lewis Caroll – Charles L. Dodgson), but writers that would do fundamental scientific or technical research don't quite exist!

# SECTION V. THE LETTERS' FREQUENCY (BY EQUAL GROUPS) IN THE ROMANIAN JURIDICAL TEXTS

Analyzing the deterioration's degree of the keys of a typing machine which functioned for more than 40 years at the clerk's office of a court of a Romanian District (Vâlcea), one partitions them in the following groups:

1) Letters completely deteriorated (one cannot read anything anymore on the typewriter).
2) Letters from which one sees only one point, hardly perceptible.
………………………..
10) Letters from which is missing only one point.
11) Letters, which are seen perfectly, without anything missing.
12) Letters which, almost have not been touched, being covered with dust.

The following resultants were obtained:

1) E, A        7) O, C, U, D, Z
2) I           8) N
3) R           9) L
4) T           10) V, M
5) S           11) F, G, B, H, X, J, K
6) P           12) W, Q, Y

This classification is a little different of that of [1], because the letters A, Ă, Â are here counted as one letter: A, The same I and Î in I, S and Ș in S, T and Ț in T.

By studying the chart of this text (from [2]), we obtain:

$$\alpha(j) = \frac{1}{23}\sum_{i=1}^{23}|\alpha(A_i)| \approx 2.348$$

thus the chart of the juridical language of current frequencies is much more larger than that of the cross words language: $\alpha(g) \approx 1.391$ and $\alpha(d_r) \approx 1.185$.

The letters P, Z and N realized the most spectacular jump:

$$\alpha(P) = 6, \ \alpha(Z) = 7, \ \alpha(N) = 8.$$

Perhaps this article surprises by its banality. But, whereas other authors spent month of calculations using computers, choosing certain books and counting the letters (!) by the computer, I have deducted this frequency of the letters in a few minutes (!), by a simple observation.

SECTION VI.  LINGUISTIC-MATHEMATICAL STATISTICS IN RECENT
ROMANIAN POETRY

**"Mathematics is logical enough to be able to detect the internal logics of poetry and crazy enough not to lag behind the poetic ineffable" (Solomon Marcus).**

The author of this article aims a statistical investigation of a recently published volume of poetry [3], which will make possible some more general conclusions on the evolution of poetry in the $XX^{th}$ century (either the literary current hermetism, surrealism or any other). Certain modifications in the structure of poetry, occurred in its evolution from classicism to modernism, are also presented. Men of letters have never agreed with mathematics and, especially, with its interference in art. Let us quote one of them: "Remarque que, a mon avis, tout literature est grotesque…(…) La seule excuse de l'écrivain c'est de se rendre compte qu'il joue, que la littérature est un jeu" (Eugène Ionesco). Well, if literature is a game why could not be subjected to mathematical investigation?

The book chosen for this study (see [3]) contains 44 poems (from which the first and the last are sort of poems essays on Romanian poetry). It comprises over 250 sentences, over 700 verses, over 2,500 words and over 11,700 letters (not sounds).

# MORPHOLOGICAL ASPECTS

1. The frequency of words depending on the grammatical category they belong to.

| | | |
|---|---|---|
| 1. Nouns | 35.592% | |
| 2. Verbs (predicat.moods) | 13.079% | "Empty" words |
| 3. Adjectives | 6.183% | 40.271% |
| 4. Adverbs | 4.829% | |
| "Full" words | 59.729% | |

1. The "full" words category includes – according to the author – nouns, verbs (predicative moods only), adjectives and adverbs. The "empty" words category includes verbs (i,e, infinitives, gerunds, poet participles, supines), numerals, articles, pronouns, conjunctions, prepositions and interjections. The same terminology was also used by Solomon Marcus in his "Poetica matematica" published by Ed. Academiei, Bucharest, 1970 (it was translated in German and published by Athenäum, Frankfurt-am-Mein, 1973).

2. The average distribution of "full" words[1] per verses (lines), sentences, poems

|   |   |
|---|---|
| a) 1.255 | nouns/line |
| b) 0.461 | verbs (p.m)/line |
| c) 0.218 | adjectives/line |
| d) 0.172 | adverbs/line |
| e) 3.464 | nouns/sentence |
| f) 1.273 | verbs (p.m)/sentence |
| g) 0.602 | adjectives/sentence |
| h) 0.475 | adverbs/sentence |
| i) 20.393 | nouns/poem |
| j) 7.492 | verbs (p.m)/poem |
| k) 3.543 | adjectives/poem |
| l) 2.792 | adverbs/poem |

We may conclude:

CONJECTURE 1. In the recent Romanian poetry the percentage of adjectives is, on average, under that of the total of words.

CONJECTURE 2. The percentage of verbs (predicative moods) is., on average, under 15% of the total of the total words.

In support of conjectures 1 and 2 we also mention:

- only one in six nouns is modified by an adjective, i.e. the role of the adjective diminishes and there are poems with no adjectives (see [3], pp. 9, 12, 20);

- on average, there is one verb in a predicative mood in more than two lines, i.e. the role of the verbal predicate decreases and there are poems with no verbal predicates (see [3], p. 20);

(From classicism to modernism both adjectives and verbal predicates gradually but constantly regressed).

- the poetry of the young poets is characterized by economy of words and, implicitly, by the avoidance of the overused words; the adjectives were favored by the romantics and the young poets feel the necessity to "renew" poetry;

- this renewal and effort to avoid the trivial may be also helped by elimination of adjectives. The strict use of adjectives or verbal predicates is also accounted for by the characteristics of the two main literary currents of our century.

a) hermetism – appeared after World War I – consists, mainly in the hyper intellectualization of language and its codification; an adjective (i.e. an explanation concerning an object) or the predicative mood of a verb (strict definition of the grammatical tense) may diminish the degree of ambiguity, generalization or abstraction intended by the poet.

b) Surrealism – literary of vanguard – aimed at detecting the irrational, the unconscious, the dream; because of its precise definite character, the adjective makes the reader "plunge" into the so carefully avoided real world.

CONJECTURE 3. In the recent Romanian poetry percentage of "full" words is over 55% of the total words.

Unlike in the spoken language in which the percentage of "full" and "empty" words is equal (see [1]) in poetry the percentage of "full" words is greater. This is due to the fact that poetry is essence, it is dense, concentrated. The percentage of "full" words and the "density" of a literary work are directly proportional.

As a conclusion to the three conjectures we may say that:

- in its evolution from classicism to modernism the percentage of nouns increased, while that of verbs decreased, less adverbs are used, on the other hand, because of the smaller number of verbs. In all, however, the percentage of "full" words increased.

3. The frequency of the nouns with and without an article.

| | |
|---|---|
| 1. Percentage of nouns with an article | - 47.884% |
| 2. Percentage of nouns without an article | - 52.116% |

CONJECTURE 4. In the recent Romanian poetry the number of nouns with an article is, on an average, smaller than the number of those without an article. With an article the noun is more definite, specified which are characteristics undesirable from the same viewpoint as that mentioned above. That is why the indefinite article is favored in modern poetry. The consequence of this preferred indefinite character of the noun enlarges the abstraction, generalization, ambiguity and, hence, the "density" of the poem. (See also the second part of assertions 1 and 2 and the statistical conjecture 3). In its evolution from classicism to modernism the number of nouns without an article used in poetry also increased.

4. The frequency of nouns depending on the grammatical case they belong to.

| Nominative | Genitive | Dative | Accusative | Vocative |
|---|---|---|---|---|
| 29.497% | 19.888% | 0.335% | 50.056% | 0.224% |
| 2 | 3 | 4 | 1 | 5 |
| ↑C L A S S I F I C A T I O N↑ | | | | |

CONJECTURE 5. In the poems under study, over 75% of the nouns are accusative or nominative.

5. Sentences, lines, words, syllables, letters – average relationships

|  |  |
|---|---|
| a) 2.402 | letters/syllable |
| b) 1.933 | syllables/word |
| c) 4.643 | letters/word |
| d) 3.528 | words/line |
| e) 6.820 | syllables/line |
| f) 16.380 | letters/line |
| g) 2.760 | lines/sentence |
| h) 9.737 | words/sentence |
| i) 18.823 | syllables/sentence |
| j) 45.208 | letters/sentence |
| k) 5.887 | sentences/poem |
| l) 16.250 | lines/poem |
| m) 57.330 | words/poem |
| n) 110.825 | syllables/poem |
| o) 266.175 | letters/poem |

Conclusion: the poems are of medium length; the lines are short while the sentences are, again, of medium length.

6. The frequency of words according to their length (in syllables)

| Syllables | Percentages | Order |
|---|---|---|
| 1 | 41.509% | 1 |
| 2 | 32.069% | 2 |
| 3 | 19.363% | 3 |
| 4 | 5.688% | 4 |
| 5 | 1.371% | 5 |
| 6 | 0.000% | 6 |

The total number of syllables in the volume is … 4,800. The frequency of words and their length (in syllables are in inverse ratio. Long words seem "less poetical".
CONJECTURE 6. In the recent Romanian poetry the percentage of words of one and two syllables is … 75%. Again, it seems that short and very short words (of one and two syllables) appear more adequate to satisfy the internal rhythm of the poem. Longer words already have their own rhythm dictated by the juxtaposition of the syllables; it is very probable that this rhythm comes into … with the rhythm imposed by the poem. Shorter words are more easily uttered; longer words seem to render the text more difficult.
7. The frequency of words according to their length (in letters)

| 1 letter | 2 | 3 | 4 | 5 | 6 | 7 | 8 | 9 | 10 | 11 | 12 | 13 | 14 |
|---|---|---|---|---|---|---|---|---|---|---|---|---|---|
| 3.604% | 25.426% | 8.475% | 11.089% | 13.347% | 13.149% | 13.703% | 5.861% | 3.129% | 1.149% | 0.752% | 0.237% | 0.079% | 0.000% |
| Order 8 | 1 | 6 | 5 | 3 | 4 | 2 | 7 | 9 | 10 | 11 | 12 | 13 | 14 |

In the whole volume there are only two words of 13 letters and 6 of twelve. A 90% of the words consist of no more than 7 letters.

CONJECTURE 7. In the recent Romanian poetry the percentage of the two letter words is, on average, about 25% of the words. In fact, the same percentage, or even higher, is found in the ordinary language. Because of esthetic reasons in poetry there is a slight tendency of reducing the frequency of the two letter words – which are especially, prepositions and conjunctions.

8. The frequency of the letters:

| The order of the letter | Letter | The average % of the frequency of the letter | The average % of vowels | The average % of cons |
|---|---|---|---|---|
| 1 | E | 11.994% | | |
| 2 | I | 10.166% | | |
| 3 | A | 8.406% | | |
| 4 | R | 7.680% | | |
| 5 | N | 6.407% | | |
| 6 | U | 6.347% | | |
| 7 | T | 5.792% | | |
| 8 | L | 5.237% | | |
| 9 | C | 5.143% | 46.865% | |
| 10 | S | 4.220% | | |
| 11 | O | 3.699% | | |
| 12 | P | 3.451% | | |
| 13 | Ă | 3.417% | | 53.135% |
| 14 | M | 3.178% | | |
| 15 | D | 2.981% | | |
| 16 | Î | 2.828% | | |
| 17 | V | 1.435% | | |
| 18 | G | 1.48% | | |
| 19 | B | 1.358% | | |
| 20 | Ș | 1.281% | | |
| 21 | F | 1.179% | | |
| 22 | Z | 0.846% | | |
| 23 | Ț | 0.803% | | |
| 24 | H | 0.496% | | |
| 25 | J | 0.196% | | |

| 26 | X | 0.034% | | |
|---|---|---|---|---|
| 27 | Ă | 0.008% | | |
| 28-31 | K | 0.000% | | |
| 28-31 | Q | 0.000% | | |
| 28-31 | Y | 0.000% | | |
| 28-31 | W | 0.000% | | |

CONJECTURE 8. In the recent Romanian poetry the percentage of vowels is, on average, over 45% of the total of letters.

Explanation: in the ordinary language the percentage of vowels is 42.7% (see [1]). In poetry it is greater because:

- vowels are more "musical" than consonants; therefore the words with more vowels "seem" more poetical; words with many vowels confer a special sonority to the text;

- modern poets and poetry are more preoccupied by form than by content, so that more attention is given to expression; the form may prejudice the content, because, very often, the reader is "caught" by sonority and less by essence;

- the internal rhythm of poetry, usually absent in the ordinary language, is also conditioned, partially, by a greater number of vowels;

- rhyme, when used, also favors a greater percentage of vowels. The percentage of vowels was greater in the period of classicism of poetry when the rhythm and rhyme were more frequently used. The special requirements of poetry impose a thorough filtration of the ordinary language.

Given the frequency of the letters in the Romanian language [1] in general:

1. E   5. N   9. L    13. D   17. S   21. F   25. J
2. I   6. T   10. S   14. P   18. B   22. T   26. X
3. A   7. T   11. O   15. M   19. V   23. Z   27. K
4. R   8. C   12. A   16. I   20. G   24. H

we may calculate the deviation of this volume of verses from the ordinary language:

$$\alpha(v) = \frac{1}{27} \sum_{i=1}^{27} |\alpha(A_i)| \approx 0.741$$

where $\alpha(A_i)$ is the deviation of the letter $A_i$, $1 \leq i \leq 27$.

The informational energy, according to O. Onicescu, is

$$\mathcal{E}(v) = \sum_{i=1}^{27} p_i^2 \approx 0.064,$$

where $p_i$, $1 \leq i \leq 27$, is the probability that the letter $p_i$ may appear in the volume (see [1]).

The first order entropy of the volume (according to Shannon) is:

$$H_1(v) = -\frac{1}{\log_{10} 2} \cdot \sum_{i=1}^{27} p_i \log_{10} p_i \approx 4.222.$$

9. The themes of the volume are studied by determining the recurrent elements, those that seem to obsess the poet. We will call these elements "key-words" and they are, in order: nouns, verbs, adjectives. Their frequency in the volume is studied. The more frequent words are all included in common notional spheres that will "decode" the themes dealt with by the poet in the volume under study, i.e.:

Elements of the Nature            Cosmological Elements

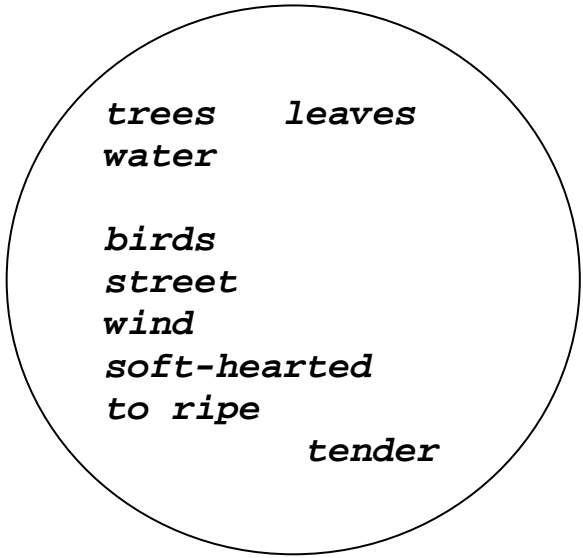
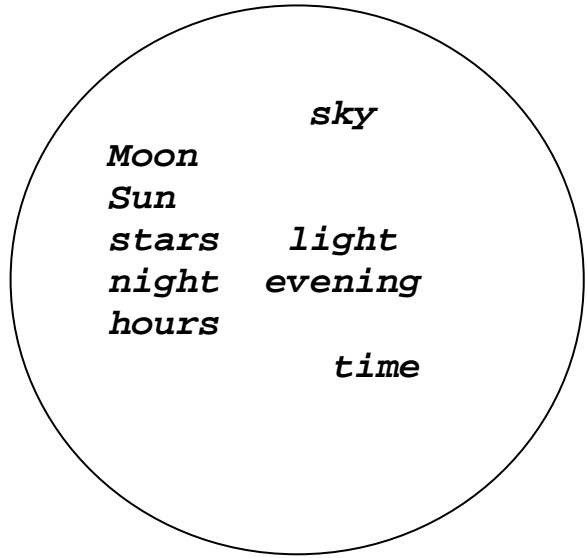

*Existence Elements*              *Poet's condition*

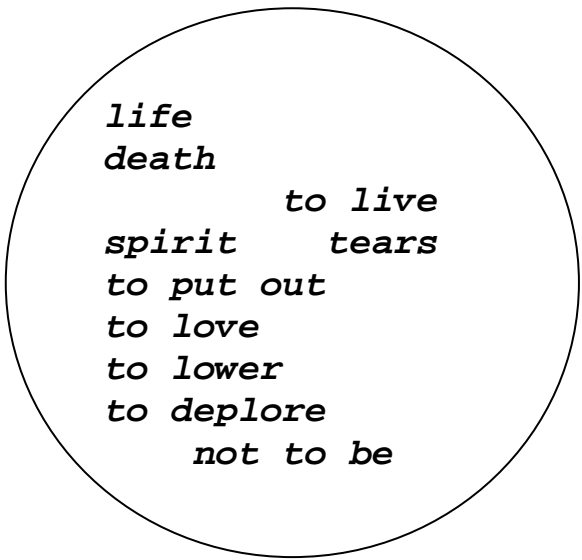
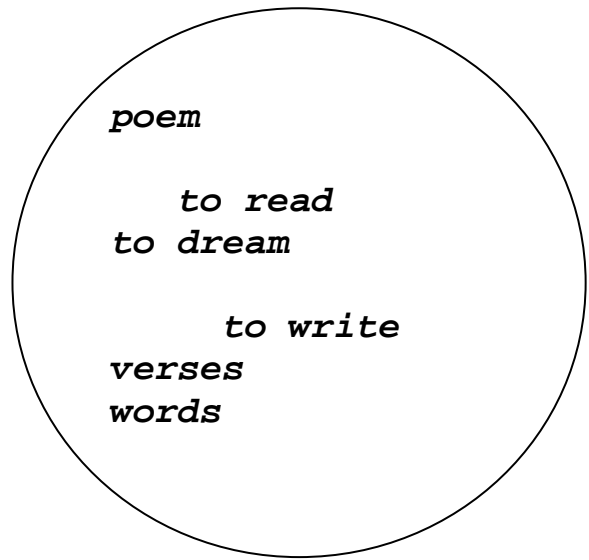

These 33 key-words (together with their synonyms) confer certain pastoral note (this was noticed by Constantin Matei, the newspaper "Înainte", Craiova), cosmological (Constantin M. Popa), existentialist nuances (Aureliu Goci, "Luceafarul", Bucharest); the preoccupation of the poet for the condition of the poet and society (Ion Pachia Tatomirescu, Craiova) is also revealed by the frequent use of certain suggestive words.

Of all the words, 33 key-words together with their synonyms have the greatest frequency in the volume.

10. The frequency of words and phrases strongly deviated from the "normal", i.e. the rules of the literary language are about 1.980 of the total of words. (We mean expressions like: "state of self", "very near myself", "it is raining at plus infinite" or words like "nontime", etc. (see [3], pp. 9, 29, 40, 31).

CONJECTURE 9. In the recent Romanian poetry the percentage of words and phrases that strongly deviated from the "normal" of the ordinary language, as well as the rules of the literary language, is slightly over 1. This fact may be accounted for by:
- content seems less important; poets are more concerned with form;
- poets invent words and expressions to be able to better reveal their feelings and emotions;
- the association of antonyms may give birth to constructions that, somehow "violate" the normal;
- poetry is, in fact, destined to break the rules and rebel against the ordinary fact (if, this right is denied, any newspaper article could be called poetry).

"In art" said Voltaire, "rules are only meant to be broken".

In its evaluation from classicism to modernism the percentage of such abnormal words and constructions increased, starting, in fact from zero. Modern literary currents favor the appearance of them.

[Editions Scientifiques, Casablanca, Morocco, 1984]

SECTION VII.  MATHEMATICAL FANCIES AND PARADOXES

**MISCELLANEA**

1. Archimedes' "fixed point theorem": Give me a fixed point in space, and I shall upset the Earth".

2. MATHEMATICAL LINGUISTICS[1]
       Poem by Ovidiu Florentin[2]

**Definition**
A word's sequence converges if it is found in a neighborhood of our heart..
*
The hermetic verses are linear equations.
*
**Theorem**
For any X there is no Y such that Y knows everything which X knows. And the reciprocal.
The proof is very intricate and long, and we will present it. We leave it to the readers to solve it!
**
My Law: Give me a point in space and I shall write the proposition behind it.

Final Motto
- O, MATHEMATICS, YOU, EXPRESSION OF THE ESSENTIAL IN NATURE!

---

1  Volume which includes this mathematical poem (pp. 39-41).
2  (Translated from Romanian by the author.) It is the mathematician's penname. He wrote many poetical volumes (in Romanian and French), as "Legi de compoziție internă. Poeme cu…probleme!" (Laws of internal composition. Poems with…problems!), Ed. El Kitab, Fès, Morocco, 1982.

## AMUSING PROBLEMS

1. Calculate the volume of a square.
(Solution: Volume = Area of the Base x Height = Side$^2$ x 0 = 0! We look at the square as an extreme case of parallelepiped with the height null.)

2. ? x 7 = 2?
(Solution: of course $\frac{2}{7} \times 7 = 2$ !)

3. Ten birds are on a fence. A hunter shoots three of them. How many birds remain?
(Answer: **none**, because the three dead birds fell down from the fence and the other seven flew away!)

4. Ten birds are in a meadow. A hunter shoots three of them. How many birds remain?
(Answer: **three birds**, the dead birds, because the others flew away!)

5. Ten birds are in a cage. A hunter shoots three of them. How many birds remain?
(Answer: **ten birds**, dead and alive, because none could get out!)

6. Ten birds are up in the sky. A hunter shoots three of them. How many birds remain?
(Answer: **seven birds**, at last, those who are still flying and those that fell down!)

7. Prove that the equation $X = X + 2$ has two distinct solutions.
(Answer: $X = \pm\infty$ !)

8. (Solving Fermat's last theorem) Prove that for any non-null integer n, the equation $X^n + Y^n = Z^n$, $XYZ \neq$, has at least one integer solution!
(Answer: (a) $n \geq 1$. Let $X_k = Y_k = Z_k = 2^k$, $k = 1, 2, 3, \ldots$ All $X_k \in N$, $K \geq 1$. $L = \lim_{k \to \infty} X_{k \in N}$. But $L = \infty \in N$, that is the integer infinite, and $\infty^n + \infty^n = \infty^n$ ! If $n$ is even, the equation has eight distinct integer solutions: $X = Y = Z = \pm\infty$ ! Similarly, we take the negative integer infinite: $-\infty \in Z-$ ]

(b) $n \leq -1$. Clearly there are at last eight distinct integer solutions: $X = Y = Z = \pm\infty$ !)

## WHERE IS THE ERROR IN THE BELOW DIOPHANTINE EQUATIONS ?

Statement:
(1) To solve in $\mathbb{Z}$ the equation: $14x + 26y = -20$.
"Resolution": The integer general solution is:
$$\begin{cases} x = -26k + 6 \\ y = 14k - 4 \end{cases} (k \in \mathbb{Z})$$

(2) To solve in $\mathbb{Z}$ the equation: $15x - 37y + 12z = 0$.
"Resolution" The integer general solution is:

$$\begin{cases} x = k+4 \\ y = 15k \\ z = 45k - 5 \end{cases} \quad (k \in \mathbb{Z})$$

(3) To solve in $\mathbb{Z}$ the equation: $3x - 6y + 5z - 10w = 0$.

"Resolution" the equation is written: $3(x - 2y) + 5z - 10w = 0$.

Since $x, y, z, w$ are integer variables, it results that 3 divides $z$ and that 3 divides $w$. I. e: $z = 3t_1$ $(t_1 \in \mathbb{Z})$ and $w = 3t_2$ $(t_2 \in \mathbb{Z})$.

Thus $3(x - 2y) + 3(5t_1 - 10t_2) = 0$ where $x - 2y + 5t_1 - 10t_2 = 0$.

Then: $\begin{cases} x = 2k_1 + 5k_2 - 10k_3 \\ y = k_1 \\ z = \quad\quad 3k_2 \\ w = \quad\quad\quad\quad 3k_3 \end{cases}$ with $(k_1, k_2, k_3 \in \mathbb{Z}^3)$,

constitute the integer general solution of the equation.

Find the error of each "resolution".

**SOLUTIONS.**

(1) $x = -26k + 6$ and $y = 14k - 4$ ($k \in \mathbb{Z}$) is an integer solution for the equation (because it verifies it), but it is not the general solution, because $x = -7$ and $y = 3$ verify the equation, they are a particular integer solution, but:

$\begin{cases} -26k + 6 = -7 \\ 14k - 4 = 3 \end{cases}$ implies that $k = \dfrac{1}{2}$ (does not belong to $\mathbb{Z}$).

Thus one cannot obtain this particular from the previous general solution.

The true general solution is: $\begin{cases} x = -13k + 6 \\ y = 7k - 4 \end{cases}$ ($k \in \mathbb{Z}$). (from [1])

(2) In the same way, $x = 5$, $y = 3$, $z = 3$ is a particular solution of the equation, but which cannot be obtained from the "general solution" because:

$\begin{cases} k + 4 = 5 \Rightarrow k = -1 \\ 15k = 3 \Rightarrow k = \dfrac{1}{5} \\ 45k - 5 = 3 \Rightarrow k = \dfrac{8}{45} \end{cases}$,

contradictions.

The integer general solution is: $\begin{cases} x = k_1 \\ y = 3k_1 + 12k_2 \\ z = 8k_1 + 37k_2 \end{cases}$ (with $(k_1, k_2) \in \mathbb{Z}^2$, cf. [1]).

(3) The error is that: "3 divides ($5z - 10w$)" does not imply that "3 divides $z$ and 3 divides $w$". If one believes that one loses solutions, then this is true because $(x, y, z, w) = (-5, 0, 5, 1)$ constitutes a particular integer solution, which cannot be obtained from the "solution" of the statement.

The correct resolution is: $3(x - 2y) + 5(z - 2w) = 0$, that is $3p_1 + 5p_2 = 0$, with $p_1 = x - 2y$ in $\mathbb{Z}$, and $p_2 = z - 2w$ in $\mathbb{Z}$.

It results that: $\begin{cases} p_1 = -5k = x - 2y \\ p_2 = 3k = z - 2w \end{cases}$ in $\mathbb{Z}$.

From which one obtains the integer general solution:
$$\begin{cases} x = 2k_1 - 5k_2 \\ y = k_1 \\ z = 3k_2 + 2k_3 \\ w = k_3 \end{cases} \text{ with } (k_1, k_2, k_3) \in \mathbb{Z}^3$$

[1] One can find these solutions using: Florentin SMARANDACHE - "Un algorithme de résolution dans l'ensemble des numbers entiers pour les équations linéaires".

## WHERE IS THE ERROR IN THE BELOW INTEGRALS ?

Let the function $f : \mathbb{R} \to \mathbb{R}$ be defined by $f(x) = 2 \sin x \cos x$.
Let us calculate its primitive:
(1) First method.
$$\int 2 \sin x \cos x \, dx = 2 \int u \, du = 2 \frac{u^2}{2} = u^2 = \sin^2 x, \text{ with } u = \sin x.$$
One thus has $F_1(x) = \sin^2 x$.
(2) Second method:
$$\int 2 \sin x \cos x \, dx = -2 \int \cos x (-\sin x) dx = -2 \int v \, dv = -v^2,$$
thus $F_2(x) = -\cos^2 x$
(3) Third method:
$$\int 2 \sin x \cos x \, dx = \int \sin 2x \, dx = \frac{1}{2} \int (\sin 2x) \, 2dx = \frac{1}{2} \int \sin t \, dt = -\frac{1}{2} \cos t$$
thus $F_3(x) = -\frac{1}{2} \cos 2x$.

One thus obtained 3 different primitives of the same function.
How is this possible?
*Answer*: There is no error! It is known that a function admits an infinity of primitives (if it admits one), which differ only by one constant.
In our example we have:
$F_2(x) = F_1(x) - 1$ for any real $x$, and $F_3(x) = F_1(x) - \frac{1}{2}$ for any real $x$.

## WHERE IS THE ERROR IN THE BELOW REASONING BY RECURRENCE ?

At an admission contest at an University, was given the following problem:
"Find the polynomials $P(x)$ with real coefficients such that $xP(x-1) = (x-3)P(x)$, for all $x$ real."

Some candidates believed that they would be able to show by recurrence that the polynomials of the statement are those which verify the following property: $P(x) = 0$ for all natural values.

In fact, they said, if one puts $x = 0$ in this relation, it results that $0 \cdot P(-1) = -3 \cdot P(0)$, therefore $P(0) = 0$.

Likewise, with $x = 1$, one has: $1 \cdot P(0) = -2 \cdot P(1)$, therefore $P(1) = 0$, etc.

Let's suppose that the property is true for $(n-1)$, therefore $P(n-1) = 0$, and we are looking to prove it for $n$:

One has: $n \cdot P(n-1) = (n-3) \cdot P(n)$, and since $P(n-1) = 0$, it results that $P(n) = 0$.

Where the proof failed?

*Answer:* If the candidates would have checked for the rank $n = 3$, they would have found that: $3 \cdot P(2) = 0 \cdot P(3)$ thus $0 = 0 \cdot P(3)$, which does not imply that $P(3)$, is null: in fact this equality is true for any real $P(3)$.

The error, therefore, is created by the fact that the implication: "$(n-3) \cdot P(n) = n \cdot P(n-1) = 0 \Rightarrow P(n) = 0$" is not true.

One can find easily that $P(x) = x(x-1)(x-2)k$, $k \in \mathbb{R}$.

## WHERE IS THE ERROR?

Given the functions $f, g : \mathbb{R} \to \mathbb{R}$, defined as follows:

$$f(x) = \begin{cases} e^x, & x \leq 3 \\ e^{-x}, & x > 3 \end{cases} \quad \text{and} \quad g(x) = \begin{cases} x^2, & x \leq 0 \\ -2x + 7, & x > 3 \end{cases}$$

Compute $f \circ g$.

"*Solution*": We can write:

$$f(x) = \begin{cases} e^x, & x \leq 0 \\ e^x, & 0 < x \leq 3 \\ e^{-x}, & x > 0 \end{cases} \quad \text{and} \quad g(x) = \begin{cases} x^2, & x \leq 0 \\ -2x + 7, & 0 < x \leq 3 \\ -2x + 7, & x > 3 \end{cases}$$

from where

$$(f \circ g)(x) = f(g(x)) = \begin{cases} e^{x^2}, & x \le 0 \\ e^{-2x+7}, & 0 < x \le 3 \\ e^{2x-7}, & x > 3 \end{cases}$$

and $f \circ g : \mathbb{R} \to \mathbb{R}$.

Correct solution:

$$f \circ g = f(g(x)) = \begin{cases} e^{g(x)}, & \text{if } g(x) \le 3 \\ e^{-g(x)}, & \text{if } g(x) > 3 \end{cases} \qquad f \circ g : \mathbb{R} \to \mathbb{R}$$

$$g(x) \le 3 \Rightarrow \begin{cases} x^2 \le 3 \Rightarrow x \in [-\sqrt{3}, 0] \\ \text{or} \\ -2x + 7 \le 3 \Rightarrow x \in [2, +\infty) \end{cases}$$

$$g(x) > 3 \Rightarrow \begin{cases} x^2 > 3 \Rightarrow x \in (-\infty, -\sqrt{3}) \\ \text{or} \\ -2x + 7 > 3 \Rightarrow x \in (0, 2) \end{cases}$$

Therefore

$$f \circ g)(x) = \begin{cases} e^{-x^2}, & x \in (-\infty, -\sqrt{3}) \\ e^{x^2}, & x \in [-\sqrt{3}, 0) \\ e^{2x-7}, & x \in (0, 2) \\ e^{2x-7}, & x \in [2, +\infty) \end{cases}$$

[Published in "Gazeta matematică", nr.7/1981, Anul LXXXVI, pp.282-283.

**WHERE IS THE ERROR IN THE BELOW SYSTEM OF INEQUALITIES ?**

Solve the following inequalities system:
$$\begin{cases} x \ge 0 & (1) \\ y \ge 0 & (2) \\ x - 2y + 3z \ge 0 & (3) \\ -3x - y + 4z \ge 4 & (4) \end{cases}$$

*"Solution"*: Multiply the third inequality by 3 and add it to the fourth inequality. The sense will be conserved. It results:
$-7y + 13z \ge 4$, or $z \ge \dfrac{1}{13}(7y + 4)$.

Therefore, $x \geq 0$ and $y \geq 0$ (from the inequalities (1) and (2))

and $z \geq \frac{1}{13}(7y+4)$     (*).

But $x = 13 \geq 0$, $y = 0 \geq 0$, and $z = 2 \geq \frac{4}{13} = \frac{1}{13}(7 \cdot 0 + 4)$ verifies (*). But we observe that it does not verify the inequalities system, because substituting in the fourth inequality we obtain: $-3 \cdot 13 - 0 + 4 \cdot 2 \geq 4$ which is not true.

Where is the contradiction?

*Solution.*

The previous solution is incomplete. We didn't intersect all four inequalities. Giving a geometrical interpretation in $\mathbb{R}^3$, and writing the inequalities as equations, we have, in fact, four planes, each dividing the space in semi spaces. Therefore, the system's solution will be formed by the points which belong to the intersection of those four semi spaces, (each inequality determines a semi space). The inequality obtained by adding the third inequality with the fourth represents, is, in fact, another semi space that includes the system's solution, and it does not simplify the system (in the sense that we cannot eliminate any of the system's inequalities).

Therefore $x = 0$, $y = 3$, $z = \frac{5}{13}$ verifies (*) but it does not verify, this time, the third inequality (although the fourth one is verified).

## THE ILLOGICAL MATHEMATICS!

Find a "logic" with the following statements:
(1) $4 - 5 \approx 5$!
(2) 8 divided by two is equal to zero!
(3) 10 minus 1 equals 0.
(4) $\int f(x) \, dx = f(x)$!
(5) 8+8=8!

*Solutions:*

These mathematical fantasies are entertainments, amusing problems; they disregard current logic, but having their own "logic", fantasist logic: thus
   (1) can be explained if one does not consider "4 - 5" as the writing of "4 minus 5" but that of "from 4 to 5"; from which a reading of the statement "$4 - 5 \approx 5$" should be: "between 4 and 5, but closer to 5".
   (2) 8 can be divided by two … in the following way:…,  i. e. it will be cut into two equal parts, which are equal to "0" above and below the cutting line!
   (3) "10 minus 1" can be treated as: the two typographical characters 1, 0 minus the 1, which justifies that there remains the character 0.
   (4) The sign will be considered as the opposite function of the integral.
   (5) The operation "$\infty + \infty = \infty$" is true: writing it vertically:

$$\infty + \infty = \infty$$

which, transposed horizontally (by a mechanic rotation of the graphic signs) will give us the statement: " 8 + 8 = 8 ".

## OPTICAL ILLUSION
### (Mathematical Psychology)

What digit is it, 8 or 3?

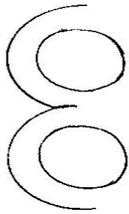

[Answer: Both of them!]

1. EPMEK = Reverse of Kempe.

2. DEDE/KIND = DedeKind's cut.

3.
```
          B
      R
    O             =    Angle of Brocard.
      C
        A
          R
            D
```

4. •BRIANCHON = Point of Brianchon.

5. $\begin{vmatrix} SYL \\ VES \\ TER \end{vmatrix}$ = Determinant of Sylvester.

6. E A O T E E
   r  t  sh n s  = The Sieve of Eratosthenes.

7.  
        A  
   R    C  
    T   S  
  D  E  S      =      Foliate curve of Descartes.

8. $\begin{pmatrix} MRX \\ RAI \\ XIT \end{pmatrix}$ = Symmetrical matrix.

9. $\overline{\text{SHEFFER}}$ = Bar of Sheffer.

10. 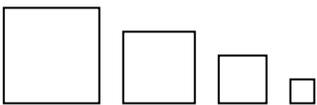 = Method of the smallest squares.

11. $\begin{pmatrix} J10000 \\ 0\varnothing1000 \\ 00R100 \\ 000D10 \\ 0000A1 \\ 00000N \end{pmatrix}$ = Matrix of Jordan.

12. NOITCNUF = Inverse function.

13. SERUGIF = Inverse figures.

14.  
     R  V  R  V  
  M  K  M  K    =     Markov Chains.  
    A  O  A  O

15. $\dfrac{\text{USA}}{\text{WEST EUROPE}}$ = Harmonious rapport.

16. $\dfrac{\text{USA}}{\text{USSR}}$ = Unharmonious rapport.

17. 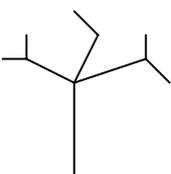 = Tree.

18. 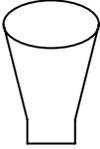 = Convergent filter.

19.
```
        A
    P       S
  O           U      = Apollonius' circle.
    L       I
        O N
```

20. 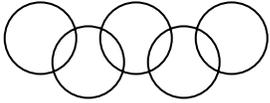 = Fascicles of circles.

21. 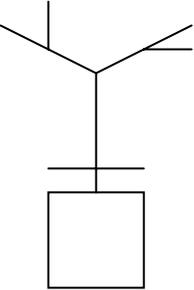 = Square root.

22. 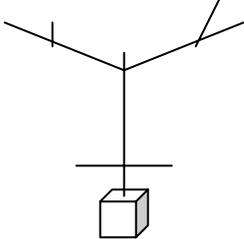 = Cubic root

23. $X^\infty + Y^\infty = Z^\infty$ = Fermat's last theorem

24. I-W-A-S-A-W-A = Iwasawa's decomposition

25.  R E
               = Latin square!
     O M

26. 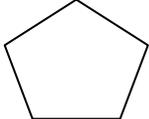 = The Pentagon!

27.  Ø                     =   Reductio ad absurdum.

28.  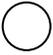   =   Ring.

29.  F           N
     U       O             =   Convex function.
        N   I
          C T

30.  P    N    S            =   Non-collinear points.
        I   T
     O

31.       G
     R      P               =   Group of rotations.
        O U

32.  ELEMENTS              =   Non-disjoint elements.

33.       M
     X       A              =   Circular matrix.
     I   T
          R

          O L
34.  P           I          =   7-gon.
     N
        O G

35.     SPA
        CE                  =   Compact space.

36.     A
        L
        G
        E                   =   Higher algebra
        B
        R
        A

37.  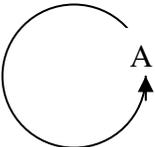   =   Vicious circle

38. A R I T H M E T I C = The higher arithmetic.

39. ∟ = Square angle.

40. SYMBOL OF (LEOPOLD) KRONECKER = L.K.

41. KOLMOGOROV'S SPACE = USSR.

42. LANGUAGE OF CHOMSKY = American.

43. GRAMMAR OF KLEENE = English.

44. CATASTROPHIC POINT = Atom bomb.

45. MACHINE OF TURING = Motor car.

46. NUMBER OF GOLD = 79 (Chemically).

47. FLY OF LA HIRE = Insect.

48. MOMENT OF INERTIA = Apathy.

49. AXIOM OF SEPARATION = Divorce.

50. CLOSED SET = Prisoners.

51. RUSSIAN MULTIPLICATION = Conquest.

52. SLIPS OF MÖBUS = Bathing trunks.

53. SINGULAR CARDINAL = Mazarin (1602-1661, France).

54. CLAN OF LEBESGUE = His family.

| | | |
|---|---|---|
| 55. SPHERE OF RIEMANN | = | Head. |
| 56. MATHEMATICAL HOPE | = | Fields prize. |
| 57. CRITICAL WAY | = | Slope. |
| 58. BOTTLE OF KLEIN | = | Beer bottle. |
| 59. CONSTANT OF EULER | = | Mathematics. |
| 60. CONTRACTING FUNCTION | = | Frost. |
| 61. BILINEAR COMBINATION | = | Concubine. |
| 62. HARDY SPACE | = | England. |
| 63. INTRODUCTION TO ALGEBRA! | = | AL. |
| 64. INTRODUCTORY ECONOMETRICS | = | ECO. |
| 65. BOREL BODY | = | Corpse. |
| 66. CHOICE FUNCTION | = | Marriage. |
| 67. GEOMETRICAL PLACES | = | ATHENA, ERLANGEN, etc. |

[Published in GAMMA, Anul IX, No. 1, November 1986]

## MATHEMATICAL LOGIC

How many propositions are true and which ones from the following:
1. There exists one false proposition amongst those n propositions.
2. There exist two false propositions amongst those n propositions.
……………………………………………………………………….
…There exist i false propositions amongst those n propositions.
……………………………………………………………………….
n. There exist n false propositions amongst those n propositions.

(This is a generalization of a problem proposed by professor FRANCICO BELLOT, in the journal NUMEROS, nr. 9/1984, p. 69, Canary Island, Spain.)

Comments:

Let $P_i$ be the proposition $i$, $1 \leq i \leq n$. If $n$ is even, then the propositions $1, 2, \ldots, \frac{n}{2}$ are true and the rest are false. (We start our reasoning from the end; $P_n$ cannot be true, therefore $P_1$ is true; then $P_{n-1}$ cannot be true, then $P_2$ is true, etc.)

*Remark*: If $n$ is odd we have a **paradox**, because if we follow the same solving method we find that $P_n$ is false, which implies that $P_1$ is true; $P_{n-1}$ false, implies that $P_2$ is true,…, $P_{\frac{n+1}{2}}$ false implies $P_{n+1-\frac{n+1}{2}}$ true, that is $P_{\frac{n+1}{2}}$ false implies $P_{\frac{n+1}{2}}$ true, which is absurd.

If $n = 1$, we obtain a variant of liar's paradox ("I lie" is true or false?)

> 1. There is a false proposition in this rectangle.

Which is obviously a paradox.

## PARADOX OF RADICAL AXES

*Property:* The radical axes of $n$ circles in the same plan, taken two by two, whose centers are not aligned, are convergent.

"Proof" by recurrence on $n \geq 3$.

For the case $n = 3$ it is known that 3 radical axes are concurrent in a point which is called the radical center. One supposes that the property is true for the values smaller or equal to a certain $n$.

To the $n$ circles one adds the $(n+1)$-th circle.

One has (1): the radical axes of first $n$ circles are concurrent in M.

Let us take 4 arbitrary circles, among which is the $(n+1)$-th.

Those have the radical axes convergent, in conformity with the recurrence hypothesis, in the point M (since the first 3 circles, which belong to $n$ circles of the recurrence hypothesis, have their radical axes concurrent in M).

Thus the radical axes of $(n+1)$ circles are convergent, which shows that the property is true for all circles $n \geq 3$ of N.

AND YET, one can build the following counterexample:

Consider the parallelogram $ABCD$ which does not have any right angle.

Then one builds 4 circles of centers $A, B, C$ and $D$ respectively, and of the same radius. Then the radical axes of the circles $e(A)$ and $e(B)$, respectively $e(C)$ and $e(D)$, are two lines, which are medians of the segments $AB$ and $CD$ respectively.

Because $(AB)$ and $(CD)$ are parallel, and that the parallelogram does not have any right angle, it results that the two radical axes are parallel, i.e. they never intersect.

Can we explain this (apparent!) contradiction with the previous property?
Response: The "property "is true only for $n = 3$. However in the demonstration suggested one utilizes the premise (distorted) according to which for $m + 4$ the property would be true. To complete the proof by recurrence it would have been necessary to be able to prove that $P(3) \Rightarrow P(4)$, which is not possible since $P(3)$ is true but the counterexample proves that $P(4)$ is false.

## A CLASS OF PARADOXES

Let A be an attribute and non-A its negation.

P1. ALL IS "A", THE "NON-A" TOO.
Examples:
$E_{11}$: All is possible, the impossible too.
$E_{12}$: All are present, the absentee too.
$E_{13}$: All is finite, the infinite too.

P2. ALL IS "NON-A", THE "A" TOO.
Examples:
$E_{21}$: All is impossible, the possible too.
$E_{22}$: All are absent, the present too.
$E_{23}$: All is infinite, the finite too.

P3. NOTHING IS "A" NOT EVEN THE "A".
Examples:
$E_{31}$: Nothing is perfect, not even the perfect.
$E_{32}$: Nothing is absolute, not even the absolute.
$E_{33}$: Nothing is finite, not even the finite.

Remark: $P1 \Leftrightarrow P2 \Leftrightarrow P3$.
More generally: ALL (<u>verb</u>) "A", the "NON-A" too.
Of course, from these appear unsuccessful paradoxes, but the proposed method obtains beautiful ones.
Look at a pun, which reminds you of Einstein:
All is relative, the (theory of) relativity too! So:
The shortest way between two pints is the meandering way!
The unexplainable is, however, explained by this word: "unexplainable"!